\documentclass[twoside]{article}
\usepackage{eucal ,amsmath}
\usepackage{amsthm}
\usepackage{amssymb}
\usepackage{amsfonts}
\usepackage[misc]{ifsym}
\usepackage{setspace,caption}
\usepackage{setspace}
\begin{document}
\thispagestyle{empty}

\title{Order Cancellation Law in a Semigroup of Closed Convex Sets}%[Order Cancellation Law in a Semigroup of Closed Convex Sets]
\author{Jerzy Grzybowski\footnote{Corresponding author} \,and Hubert Przybycie\'n}
\date{}
\maketitle

\newcommand{\K}{\ensuremath{\mathcal{K}}}
\newcommand{\PP}{\ensuremath{\mathcal{P}}}
\newcommand{\RR}{\ensuremath{\mathbb{R}}}
\newcommand{\NN}{\ensuremath{\mathbb{N}}}
\newcommand{\B}{\ensuremath{\mathcal{B}}}
\newcommand{\SSS}{\widetilde{S}}
\newcommand{\st}{\stackrel{\ast}{+}}
\newcommand{\koniec}{${\square}$}
\newcommand{\dow}{${\mathbf{Proof.\;}}$}

\def\dotminus{\mathbin{\ooalign{\hss\raise1ex\hbox{.}\hss\cr
  \mathsurround=0pt$-$}}}

%\renewcommand\rightmark{[Order Cancellation Law in a Semigroup of Closed Convex Sets]}
%\renewcommand\leftmark{[running authors]}

%Jerzy Grzybowski\footnote{Corresponding author} and Hubert Przybycie\'n

\begin{abstract}
\noindent
In this paper generalize Robinson's version of an order cancellation 
law for subsets of vector spaces in which we cancel 
by unbounded sets. We introduce the notion of weakly narrow 
sets in normed spaces, study their properties and prove the 
order cancellation law where the canceled set is weakly narrow. Also 
we prove the order cancellation law for closed convex subsets of 
topological vector space where the canceled set has bounded Hausdorff-like 
distance from its recession cone. We topologically embed the 
semigroup of closed convex sets sharing a recession cone having 
bounded Hausdorff-like distance from it into a topological vector 
space. This result extends Bielawski and Tabor's generalization of R\aa dstr\"om 
theorem.
\end{abstract}

{\small 2010 \textit{Mathematics Subject classification}. 52A07, 18E20, 46A99.}

{\small \textit{Key words and phrases}.
semigroup of convex sets, order cancellation law, recession cone, Minkowski--R{\aa}dstr{\"o}m--H{\"o}rmander space.}

\newtheorem{tw}{Theorem}[section]
\newtheorem{Prop}[tw]{Proposition}
\newtheorem{Lem}[tw]{Lemma}
\newtheorem{Cor}[tw]{Corollary}
\newtheorem{Def}[tw]{Definition}
\theoremstyle{remark}
\newtheorem{example}[tw]{Example}
\newtheorem{Rem}[tw]{Remark}

\section{Introduction}
%\begin{section}{Introduction}
The order cancellation law or R\aa dstr\"om cancellation theorem \cite{mp,hR}, investigated also in \cite{GKKU,GU,KP}, 
enables an embedding of a semigroup of convex sets into a quotient vector space 
called a Minkowski--R{\aa}dstr{\"o}m--H{\"o}rmander space \cite{lD, rU}. 
This space plays a crucial role in differentiation of nonsmooth functions 
(quasidifferential calculus of Demyanov and Rubinov \cite{mA,DR1,DR2,mD,aF,KUGKATS}) and 
in integration and differentiation in the theory of multifunctions \cite{debreu1, debreu2}. 
The cancellation by unbounded sets is especially challenging. It is usually impossible. However, 
in this paper we significantly extend the range of applicability of this law. 

In order to embed a given semigroup of subsets of a vector space
into some other vector space, we need the law of cancellation: $A+B= B+C\Longrightarrow A= C$. 
This embedding is very important in set-valued analysis. Let us recall that  
the cancellation law holds true in the semigroup $\mathcal{B}(X)$ of nonempty bounded closed convex sets 
in real Hausdorff  topological vector space $X$, 
where the addition $\dot{+}$ is defined by $A\dot{+}B=$ cl$(A+B)=$ cl$\{a+b:a\in A, b\in B\}$. 
Minkowski studied the semigroup of compact convex subsets of Euclidean spaces, therefore, 
the algebraic or vector addition $A+B$ is also called the Minkowski addition. 
In infinite dimentional topological vector spaces, we need to consider $A\dot{+}B$ instead of $A+B$.    
R\aa dstr\"om \cite{hR} applied the cancellation law in order to embed the convex cone 
$\mathcal{B}(X)$, where $X$ is a normed vector space, 
into a normed vector space.   
H\"ormander \cite{lH} generalized the R\aa dstr\"om's result to the locally convex spaces and Urba\'nski \cite{rU} 
to the arbitrary topological vector spaces. 
However, the order law of cancellation which is applied in the aforementioned cased is based on the following theorem.  

\begin{tw}{\rm (see for example Proposition 2.1 in \cite{rU}.) }
Let $X$ be a Hausdorff topological vector space and $A,B,C\subset X$. If $B$ is nonempty and bounded and 
$C$ is closed and convex then $\textup{cl}(A+B)\subset \textup{cl}(B+C)\Longrightarrow A\subset C$.
\end{tw}

In \cite{GU} the assumption of closedness of $C$ is weakend. In this paper we drop the assumption of 
boundedness of $B$. 
 
Since the cancellation law does not hold in the semigroup $\mathcal{C}(X)$ of nonempty closed convex sets, 
this semigroup cannot be embedded into 
a vector space. However, Robinson \cite{sR} proved the cancellation law for the family $\mathcal{C}_V(\mathbb{R}^n)$ 
of nonempty closed convex 
subsets of $\mathbb{R}^n$ sharing common recession cone $V$. He also embedded $\mathcal{C}_V(\mathbb{R}^n)$ into 
a topological vector space 
(but not into a normed space). Bielawski and Tabor \cite{BT} restricted the family $\mathcal{C}_V(X)$, where $X$ is a normed vector space 
to such closed convex sets $A$ that the Hausdorff distance $d_H(A,V)$ is finite. It enabled them to embed restricted $\mathcal{C}_V(X)$ 
into a normed vector 
space, generalizing in this case the R\aa dstr\"om's result \cite{hR}.     

In this paper, we investigate the possibility of cancellation ($A+B\subset B+C\Longrightarrow A\subset C$) for some 
unbounded sets $B$. 
In particular, we generalize in Theorem \ref{rob01} and \ref{rob2} the following Robinson's result. 

\begin{tw}{\rm (see Lemma 1 in \cite{sR}.) }
Let $X=\mathbb{R}^n$ and $A,B,C\subset X$. If $B,C$ are nonempty closed and convex and 
${\rm recc}\,B ={\rm recc}\,C$ then $A+B\subset B+C\Longrightarrow A\subset C.$
\end{tw}

In Theorem \ref{gTB0} we generalize  the following Bielawski-Tabor's order law of cancellation.  

\begin{tw}{\rm  (see Proposition 1 in \cite{BT}.)} Let $X$ be a Banach space, $V$ a closed convex cone in $X$, 
$d_H$ a Hausdorff distance and 
$C_V = \{A \in \mathcal{C}(X) : d_H(A, V ) < \infty\}.$
If $A,B,C\in C_V$ then $\textup{cl}(A+B)\subset \textup{cl}(B+C)\Longrightarrow A\subset C.$
\end{tw}

Finally, we generalize topologically in Theorem \ref{BT2} Bielawski-Tabor's embedding theorem. 

\begin{tw}{\rm  (see Theorem 2 in \cite{BT}.)} Let $X$ be a Banach space, $V$ a closed convex cone in $X$, 
$d_H$ a Hausdorff distance and
$C_V = \{A \in \mathcal{C}(X) : d_H(A, V ) < \infty\}.$
Then the abstract convex cone $C_V$ can be embedded isometrically and isomorphically as a closed
convex cone into a Banach space.
\end{tw}

In Section 2 we prove an order cancellation law in an ordered semigroup with an operator of convergence. This result 
is applied in Section 5 for convex sets in topological vector spaces having certain kind of "finite distance" between a set  
and its recession cone.

In Section 3 we give some properties of asymptotic and recession cones in infinite dimensional spaces. 
We also present examples of unbounded convex sets with trivial recession cones and of a recession cone of sum of sets 
greater than a sum of recession cones of summands. 

In Section 4 we introduce the notion of narrow unbounded sets. We prove in Theorem \ref{rob01} that in the inclusion 
cl\,$(A+B)\subset$\,cl$(B+C)$ we can cancel a narrow set $B$ with an asymptotic cone contained in a recession cone 
of closed convex set $C$. We also present properties of addition of narrow sets that help us in Theorem \ref{mrh1} 
to show that the family of all weakly narrow closed convex subsets of normed vector space sharing 
a given pointed recession cone is closed with respect to Minkowski addition.

\section{Cancellation law in an ordered semigroup with an operator of convergence}
In this section, we prove the order cancellation law for, respectively, bounded and convex elements in ordered semigroups. 
The results of this section will be applied in Section 5 where we give a generalization of Bielawski-Tabor result \cite{BT} 
to topological vector spaces. 

Let $\mathbb{Q}_2^{+}$ be the set of all non-negative dyadic rationals and
let us consider a system $S=(S,+,\cdot\,\,, \leqslant,\, \lim )$, where '$+$' and '$ \cdot $' are binary operations from $S\times S$ and 
$\mathbb{Q}_2^{+} \times S$ into $S$ respectively, $ \leqslant $ is partial order on $S$ and $\lim$ is an limit operator turning 
$S$ into $\mathcal{L}^{\ast}$-space according to \cite{rE}, 1.7.18, p. 90. Specifically, 
(1) if $a_n=a$ then $\lim a_n =a$, (2) if $\lim a_n=a$ then $\lim a_{k_n} =a$ for every subsequence $(a_{k_n})$, 
(3) if $a_n$ does not converge to $a$ than there exists subsequence $(a_{k_n})$ such that no subsequence of $(a_{k_n})$ 
converges to $a$.
Moreover, let the following conditions be satisfied:
\begin{description}
\item[{\rm (S1)}] $a+(b+c) =(a+b)+c$, for all $a,b,c\in S$,
\item[{\rm (S2)}] $a+b=b+a$, for all $a,b\in S$,
\item[{\rm (S3)}] there exists an element $0\in S$ such that $a+0=a$, for all $a\in S$,
\item[{\rm (S4)}] If $a\leqslant b$ then $a+c\leqslant b+c $ for all  $a,b,c\in S$,
\item[{\rm (S5)}] $1\cdot a =a$, for all $a\in S$,
\item[{\rm (S6)}] if $a\leqslant b$ then $\alpha\cdot a\leqslant \alpha\cdot b $, for all  $a,b,\in S, \alpha\in\mathbb{Q}_2^{+}$,
\item[{\rm (S7)}] $\alpha (\beta\cdot a )=(\alpha\beta )\cdot a$, for all $\alpha , \beta \in \mathbb{Q}_2^{+} , a\in S$,
\item[{\rm (S8)}] $\alpha \cdot (a+b) =\alpha \cdot a +\alpha \cdot b$, for all $\alpha \in \mathbb{Q}_2^{+} , a,b\in S$,
\item[{\rm (S9)}] $(\alpha +\beta )\cdot a \leqslant \alpha \cdot a +\beta \cdot a$, for all  $\alpha , \beta \in \mathbb{Q}_2^{+} , a\in S$,
\item[{\rm (S10)}] for all sequnces $(a_n ), (b_n ) \subset S$ if $\lim_{n\to\infty } a_n =a \wedge \lim_{n\to\infty } b_n =b $ and $ \forall_{n\in\mathbb{N}} a_n \leqslant b_n $ then $a\leqslant b$,
\item[{\rm (S11)}]  for all sequnces $(a_n )\subset S$ and $b\in S$ if $\lim_{n\to\infty } a_n =a$ then\\ $\lim_{n\to\infty }( a_n +b ) =a+b.$
\end{description}

\begin{Def}
{\rm We say that an element $a\in S$ is}  bounded {\rm if\\ $\lim_{n\to\infty } (2^{-n}\cdot a) =0 .$}
\end{Def}
\begin{Def}\label{convexity}
{\rm We say that an element $a\in S$ is}  convex {\rm if  $(\alpha +\beta )\cdot a = \alpha \cdot a +\beta \cdot a$ for all  
$\alpha , \beta \in \mathbb{Q}_2^{+} $}.
\end{Def}
\begin{Lem}\label{olc<<}
Let $S$ be a semigroup satisfying conditions {\rm (S1--11)}. 
If $a,b,c\in S,$ $\lim_{n\to\infty} (2^{-n} b ) $ exists , $c$ is convex  and $a+b\leqslant b+c, $
then $a+\lim_{n\to\infty} (2^{-n} b) \leqslant c+\lim_{n\to\infty} (2^{-n} b ) .$
\end{Lem}
\begin{proof}
For $k\in\mathbb{N}$ we have
\begin{equation*}
\begin{split}
 2^k \cdot a+ b & \leqslant \underbrace{a+a+...+a}_{2^k - {\rm times}} +b\leqslant \underbrace{a+a+...+a}_{2^k -1 - {\rm times}} +b+c\leqslant \\& \leqslant \underbrace{a+a+...+a}_{2^k -2 {\rm times}} +b+c+c\leqslant ...\leqslant b+   \underbrace{c+c+...+c}_{2^k - {\rm times}} = \\&=b+2^k \cdot c.
\end{split}
\end{equation*}
Last equality follows from the convexity of an element $c$.
Multiplying both sides of the inequality by $2^{-k}$ and using (S6) and (S7) we obtain
$a+2^{-k} \cdot b \leqslant 2^{-k} \cdot b +c .$
Now by (S11)  we obtain
$a+\lim_{n\to\infty} (2^{-n} b )=\lim_{k\to\infty } (a+2^{-k}\cdot b ) $
and
$c+\lim_{n\to\infty} (2^{-n} b )=\lim_{k\to\infty } (c+2^{-k}\cdot b ) $
hence by (S10)  we get
$a+\lim_{n\to\infty} (2^{-n} b) \leqslant c+\lim_{n\to\infty} (2^{-n} b ) . $
\end{proof}
\begin{tw}\label{olc<}{{\bf (Order law of cancellation)}}
Let $S$ be a semigroup satisfying conditions {\rm (S1--11)}. If $a,b,c\in S,$ $b$ is bounded, $c$ is convex  and $a+b\leqslant b+c, $
then $a \leqslant c.$
\end{tw}
\begin{proof}
The teorem is a direct corollary from Lemma \ref{olc<<}.
\end{proof}
Let $S_{bc}\subset S$ be a subsemigroup of bounded convex elements.  
If the topology in $S_{bc}$ is generated by a uniformity in which the addition is strongly uniformly continuous 
then there exists a topological embedding 
of the semigroup $S_{bc}$ into a quotient group $\widetilde{S_{bc}}=S_{bc}\times S_{bc}/_{\sim}$, where $(a,b)\sim (c,d)\Longleftrightarrow 
a+d=b+c$
(see Proposition 1.2(i) \cite{rU} and Proposition 3.5 \cite{GPPU}). 

Theorem \ref{olc<} will be applied in the last section of this paper. The following trivial example shows a semigroup $S$ satisfying 
all conditions (S1-11) and containing only convex elements. 

\begin{example}{\rm 
Let $S=\mathbb{Q}^{+}_2\cup\{\infty\}$,  $\mathbb{Q}_2\cup\{\infty\}$, 
$\mathbb{R}_{+}\cup\{\infty\}$ or $\mathbb{R}\cup\{\infty\}$. 
Let addition be the normal addition and $S$ be ordered by natural order induced from the reals 
and $0\cdot A=\{0\}$ for all elements $A$ of $S$. 
We define a limit $\lim a_n=a$ in usual way. 
Obviously, the conditions (S1-11) are satisfied. Only the element $\infty$ is not bounded and all
elements of $S$ are convex.
}\end{example}

The next example shows a semigroup $S$ satisfying 
all conditions (S1-11) and containing unbounded and non-convex elements. 

\begin{example}{\rm 
Let $S$ be the family of all nonempty finite subsets of $\mathbb{Q}^{+}_2\cup\{\infty\}$, $\mathbb{Q}_2\cup\{\infty\}$, 
$\mathbb{R}_{+}\cup\{\infty\}$ or $\mathbb{R}\cup\{\infty\}$. Let addition be natural algebraic addition and $S$ be ordered by inclusion 
and $0\cdot A=\{0\}$ for all elements $A$ of $S$. 
We define a limit $\lim A_n=A$ if and only if 
every neighborhood of $A$ contains all but finite number of sets $A_n$
and $$A\subset \bigcap_{n_0=1}^{\infty}\textup{cl}(\bigcup_{n\geqslant n_0}A_{k_n})$$ 
for every increasing sequence $(k_n)$ of integers. 
It is easy to observe that the conditions (S1-11) are satisfied. Bounded elements of $S$ are sets not containing the element $\infty$. 
Convex elements of $S$ are all sets containing at most one number other than $\infty$.
}\end{example}

Another example is provided by a family of open domains.

\begin{example}{\rm 
Let $X$ be a real Hausdorff topological vector space and $S^{\ast}$ be the family of all open domains of $X$, 
that is of all interiors of closed subsets of $X$, and $S=S^{\ast}\cup \{\{0\}\}$.
Let addition be defined by $A\ddot{+}B=$ int cl$(A+B)$, the family $S$ be ordered by inclusion 
and $0\cdot A=\{0\}$. We define a limit $\lim A_n=A$ if and only if for any neighborhood $U$ of origin 
there exists $n_0$ such that $A_n\subset A+U$, $A\subset A_n+U$ for all $n\geqslant n_0$. 
The conditions (S1-11) are satisfied. Bounded elements of $S$ are bounded sets. 
Convex elements of $S$ are convex sets.
}\end{example}
\begin{Def}
{\rm A quartet $(S,\mathbb{R}_+,+,\cdot)$ satisfying equalities $(S1-3,5,7-8)$, $(\alpha +\beta )a=\alpha a+\beta a$ and 
$0\cdot a=0$ for all $\alpha, \beta\geqslant 0, a,b,c\in S$ is called an {\it abstract convex cone}.}
\end{Def}

%333333333333333333333333333333333333333
\section{Assymptotic and recession cones}
Our main results are placed in sections $4$ and $5$. Since considered sets in those sections are unbounded and convex, 
here we discuss some properties of the recession cone of convex set.

\begin{Def}\label{recc}
{\rm For a Hausdorff topological vector space $X$ and $A\subset X$ the set 
$A_{\infty}=\{\lim\limits_{n\rightarrow \infty}t_n x_n \,|\,  (x_n)\subset A, (t_n)\subset \mathbb{R}, 
t_n\downarrow 0  \, \}$ is called an {\it assymptotic cone} of $A$.
A set recc$A=\{x \in X\,|\, \forall a\in A, t\in \mathbb{R}_+: a+tx\in A\, \}$ is called a {\it recession cone} of $A$.}
\end{Def}

In case of ambiguity about exact topology we put $A^{\tau}_{\infty}$, $A^*_{\infty}$ or $A^{\|\cdot\|}_{\infty}$.
Obviously, a recession cone of a set is contained in an assymptotic cone of this set. 
We are going to express assymptotic and recession cones in terms of Minkowski subtraction. 

\begin{Def}
{\rm For a vector space $X$ and subsets $A,B\subset X$ the {\it Minkowski difference} of sets $A$ and $B$ is defined by} 
$A \dotminus B =\{x\in X : x+B\subset A\} .$
\end{Def}

The following obvious proposition gives us a useful formula of Minkowski difference of sets.

\begin{Prop}\label{AmBjakoprzekroj}
If $A,B\subset X $ then $A\dotminus B =\bigcap_{b\in B} (A-b) .$
\end{Prop}

Let  $X$ be a vector space and let $\mathcal{A}\subset 2^X $ be a family of subsets of $X.$
We say that $\mathcal{A} $ is {\it closed under translations } if for every $A\in \mathcal{A} $ and for every $x\in X$ we have $x+A\in  \mathcal{A} .$ 
We also say that $\mathcal{A} $ is {\it closed under intersections  } if for any subfamily $\{A_t\}_{t\in T}\subset\mathcal{A}$ the 
intersection $\bigcap_{t\in T } A_t $ belongs to $\mathcal{A} .$

The following proposition is an intermediate corollary from Proposition \ref{AmBjakoprzekroj}.

\begin{Cor}
Let  $X$ be a vector space and let $\mathcal{A}\subset 2^X $ be a family of subsets of $X.$ If  $\mathcal{A} $ 
is  closed under intersections  and translations then $\mathcal{A} $ is closed under the Minkowski difference.
\end{Cor}

It is well known that recc$A=A\dotminus A$ for any nonempty convex set $A$.

\begin{Lem}\label{recconesequence}
Let $X$ be a topological vector space and let $A\subset X$ be a closed and convex set then 
$A_{\infty}=A\dotminus A$.
\end{Lem}
\begin{proof}
Let $u\in A_{\infty} , x\in A$. By Definition \ref{recc} we have 
$u+x =\lim_{n\to\infty} (t_n x_n +x ) =\lim_{n\to\infty} (t_n x_n +(1-t_n )x )\in A$
for some $(x_n)\subset A,(t_n)\subset \mathbb{R}$ such that $t_n\downarrow 0$.
Hence $u+A\subset A .$
Therefore, $A_{\infty}\subset A\dotminus A .$

On the other hand, we have $A\dotminus A = {\rm recc}A\subset A_{\infty}$.
\end{proof}

\begin{Cor}
Let $X$ be a topological vector space. If a subset $A\subset X$ is closed and convex then also the cone $A_{\infty}$ is closed and convex.
\end{Cor}

The corollary follows from the fact that convexity and closedness are preserved under translation and intersection. 

It is well known (Theorem 8.4 in \cite{tR}) that if $X$ is finite dimensional normed space and $A\subset X$ 
is closed convex subset  then the condition $A_{\infty}=\mbox{recc}A=\{0\}$ implies that the set $A$ is bounded. 
Now, we show that this is no longer true in infinite dimensional normed space. To do this we need the following two lemmas. 

\begin{Lem}\label{existencelemma}
Let $X$ be a separable infinite dimensional normed space. Then there exists a sequence $(f_n )\subset X^* $ such that:
\begin{description}
\item[{\rm (i)}] $\| f_n \| =1$ for $n\in \mathbb{N} ,$
\item[{\rm (ii)}] for any $x\in X ,$ if $f_n (x) =0$ for every  $n\in \mathbb{N} ,$ then $x=0 .$
\end{description}
\end{Lem}
\begin{proof}
Let $A=\{x_n \,|\,n\in \mathbb{N}\}$ be a dense and countable set in $X.$ By Hahn-Banach theorem there exists 
a sequence $(f_n)\subset X^*$ of
continuous linear functionals such that $\|f_n \|=1$ and $f_n (x_n ) =\|x_n \|.$ 
Take any $x\in X $ and suppose that $f_n (x) =0 $ for every $n\in \mathbb{N} .$ Let $ (x_{n_k}) \subset A $ be 
a sequence such that $\|x-x_{n_k} \|<\dfrac{1}{k}$. We have  
$ \|x_{n_k} \| =f_{n_k}(x_{n_k} -x)\leqslant \|x-x_{n_k} \|<\dfrac{1}{k} .$ 
Hence  $x=\lim_{k\to\infty } x_{n_k} =0.$
\end{proof}
\begin{Lem}\label{unboundedset}
Let $X$ be infinite dimensional a separable normed space then there exists $A\subset X$ 
such that $A$ is closed convex unbounded and $A\dotminus A =\{0\} . $
\end{Lem}
\begin{proof}
Let $(f_n )\subset X^* $ be a sequence such as in the proof of Lemma \ref{existencelemma}. 
Let $A:=\{x\in X\,|\, |f_n (x) |\leqslant n, n\in\mathbb{N}\}.$
The set $A$ is closed and convex. Now, we show that $A$ is not bounded. 
Suppose contrary that there exists a $\sigma \in \mathbb{R} $ such that for any $x\in A$ we have $\|x\|\leqslant \sigma .$ 
Let us define a functional $\xi : X\rightarrow \mathbb{R}$, $\xi (x) :=\sup\limits_{n\in \mathbb{N}} |f_n (x)|n^{-1} .$
It is easy to see that $\xi $ is some new norm and $A=\{x\in X \,|\, \xi (x)\leqslant 1 \} $ is a unit ball in this norm.
Thus for any $x\in X $ we have $\|x\|\leqslant \sigma \xi (x) .$
Now take $m\in \mathbb{N} $ such that $m>\sigma $ and let 
$u\in\bigcap\limits_{1\leqslant j\leqslant m} \mbox{\rm ker} f_j  , u\neq 0$ then we have 
$$\|u\|\leqslant \sigma \xi (u) = \sigma\sup\limits_{n> m} |f_n (x)|n^{-1}
\leqslant \sigma m^{-1} \|u\| <  \|u\|$$
but this is impossible. Hence $A$ is not bounded.\\
Yet $A\dotminus A =\{0\}$, since $A$ is bounded in some other norm. 
\end{proof}

In the following example a set $A$ has trivial recession cone, while no linear functional is bounded on $A$. 

\begin{example}\label{xxx}
Let $c_{00}$ be a space of all sequences of real numbers with finite amount of nonzero terms. 
Let $\tau$ be any topology in $c_{00}$ in which all projections on axes of coordinates are continuous functionals.
Let $A:=$conv$\bigcup\limits_{n=1}^\infty([n,2n]^n\times \{0_{\mathbb{R}^\infty}\})$.
Let $A_{\tau}:=$cl$_{\tau}A$. Let $\sigma$ be a product topology in $c_{00}$. Then $\sigma$ is weaker than $\tau$.
Notice that $A_{\sigma }\cap (\mathbb{R}^n\times \{0_{\mathbb{R}^\infty}\})
=$conv$\bigcup\limits_{k=1}^n([k,2k]^k\times \{0_{\mathbb{R}^\infty}\})$. 
Hence $A_{\tau}=A$.
The set $A$ contains no half-line, and recc$A=\{0\}$. However, for any functional $f\in c_{00}^*$ 
the image $f(A)$ is unbounded. Hence $A$ is unbounded in any topology which separates points in $c_{00}$. 
\end{example}

\begin{tw}
Let $X$ be an infinite dimensional normed space. 
Then there exists a subset $A\subset X$ such that $A$ is closed convex unbounded and $A\dotminus A =\{0\} .$
\end{tw}
\begin{proof}
Since every infinite dimensional normed space $X$ contains some separable infinite dimensional subspace $Y$, 
our theorem follows from Lemma \ref{unboundedset}.
\end{proof}

The following theorem shows that in any infinitely dimensional Banach space the equality 
recc$A\,\dot{+}\,$recc$B\,=\,$recc$(A\,\dot{+}\,B)$ is false for some closed convex subsets $A$ and $B$.

\begin{tw}\label{reccA+B}
Let $X$ be an infinite dimensional Banach space. 
Then there exists subsets $A,B\in\mathcal{C}(X)$ with trivial reccesion cones such that 
\textup{recc}$(A\dot{+}B)$ is a half-line. 
\end{tw}
\begin{proof}
Let $X$ be an infinitely dimensional Banach space. By the Mazur theorem there exists a closed subspace $Y$ of $X$ 
having Schauder basis $(e_i), i=0,1,...$. Moreover (see \cite{jD}), there exists a constant $C\geqslant 1$ such that 
$|x_i|\leqslant C\|x\|$ for all $x=\sum\limits_{j=0}^\infty x_je_j\in Y, i=0,1,...$.  
Let $A= {\rm cl\,conv}\{ne_0+2^nne_n\,|\,n=0,1,2,\ldots\}$,
$B= {\rm cl\,conv}\{ne_0-2^nne_n\,|\,n=0,1,2,\ldots\}$.
Notice that the ray $\{te_0\,|\,t\geqslant 0\}$ is contained in $A+B$.
Hence this ray is contained in the recession cone recc$(A+B)$.

Let us assume that the set $A\dot{+}B$ contains some ray $\{tv\,|\,t\geqslant0\},v\neq 0$,
where $v\in X$, $v=\sum\limits_{n=0}^{\infty}v_ne_n$.
Let $y\in A+B$. 
Then $y=\sum\limits_{n=0}^{\infty}\alpha_n(ne_0+2^nne_n)+\sum\limits_{n=0}^{\infty}\beta_n(ne_0-2^nne_n)\textbf{}$,
where $\alpha_n,\beta_n\geqslant0, \sum\limits_{n=0}^{\infty}\alpha_n=\sum\limits_{n=0}^{\infty}\beta_n=1$ and 
the set $\{n\,|\,\max(\alpha_n,\beta_n)>0\}$ is finite.
Notice that
$$C\|tv-y\|\geqslant|tv_n-2^nn(\alpha_n-\beta_n)|\geqslant t|v_n|-2^nn$$
for all $n\geqslant 1$. If $v_n\neq 0$ then $\|tv-y\|>1$ for all
$t>\dfrac{2^nn+C}{|v_n|}$. Hence dist$(tv,A\dot{+}B)\geqslant 1$ which implies that $v_n=0$ for all $n\geqslant 1$, 
and we obtain the equality recc$\,(A\dot{+}B)=[0,\infty)e_0$.

We know that recc$\,A\,\cup\,$recc$\,B\subset \,\,$recc$\,(A\dot{+}B)$. Does $e_0$ belong to recc\,$A$? 
If yes then $e_0\in A$. Assume that 
there exists $y\in $ conv$\{ne_0+2^nne_n\,|\,n=0,1,2,\ldots\}$
such that $\|e_0-y\|<\dfrac{1}{3C}$. 
Let $y:=\sum_{n=1}^{\infty}\alpha_n\cdot(ne_0+2^nne_n)$.
Then
$$\max\big(\big|1-\sum_{n=1}^{\infty}\alpha_n n\big|,\sup\limits_{n\geqslant 1}(\alpha_n2^nn)\big)\leqslant
C\|e_0-y\|<\dfrac{1}{3}.$$
Hence we  obtain the inequalities $\alpha_nn<\dfrac{1}{3\cdot 2^n}, n\geqslant 1$. Then 
$\sum\limits_{n=1}^{\infty}\alpha_nn<\dfrac{1}{3}$.
But
$$1-\dfrac{1}{3}<\big|1-\sum\limits_{n=1}^{\infty}\alpha_nn\big|\leqslant C\|e_0-y\|<\dfrac{1}{3}.$$
Hence we obtain a contradiction. Then $e_0\not\in A\supset\,$recc$\,A$. Therefore, the recession cone of $A$ 
and, similarly, of $B$ is trivial.
\end{proof}

\begin{Cor}
Let $X$ be an infinite dimensional Banach space. 
Then there exists subsets $A,B\in\mathcal{C}(X)$ such that \textup{recc}$\,(A\dot{+}B)\neq \textup{recc}\,A\,\dot{+}\,\textup{recc}\,B$. 
\end{Cor}

\begin{Rem}
{\rm 
It is not difficult to see that in finite dimensional spaces a family of closed convex sets 
with a fixed common recession cone is closed with 
respect to Minkowski addition $A\dot{+}B=\textup{cl}(A+B)$ \cite{sR}.    
Theorem \ref{reccA+B} shows that in infinite dimensional spaces, 
even Hilbert spaces, a family of closed convex sets 
with fixed common recession cone 
need not to be closed with respect to Minkowski addition. 
In order to obtain such a family 
(with common recession cone and closed with respect to Minkowski addition), 
we need to introduce a condition 
expressed not in terms of recession cones.    
}
\end{Rem}

\begin{Rem}
{\rm 
Theorem \ref{reccA+B} shows that in any Banach space $X$ we cannot cancel by sets with trivial recession cone. 
For example let $A,B\subset X$ be sets from theorem \ref{reccA+B}. Then $A\dot{+}B\dot{+}\textup{\,recc\,}(A\dot{+}B)=A\dot{+}B$.
However, $A\dot{+}\textup{\,recc\,}(A\dot{+}B)\neq A$. Hence the condition for canceling by convex sets in Banach spaces 
cannot be expressed just in terms of recession cones.
}
\end{Rem}

\begin{tw}\label{reccA+B+}
Let $X$ be an infinite dimensional Banach space. 
Then there exists subsets $A,B\in\mathcal{C}(X)$ with trivial reccesion cones such that the sum
$A\dot{+}B$ is not closed. 
\end{tw}

If the space $X$ is not reflexive then we can always find two bounded closed convex sets, with not closed Minkowski sum. 
In fact, possibility of finding such two bounded closed convex sets is equivalent to non-reflexivity of the space $X$ \cite{hP}.  

\begin{proof}
Let $X$ be an infinitely dimensional Banach space. Again by the Mazur theorem we have a closed subspace $Y$ of $X$ 
with Schauder basis $(e_i), i=0,1,...$ and a constant $C\geqslant 1$ such that 
$|x_i|\leqslant C\|x\|$ for all $x=\sum\limits_{j=0}^\infty x_je_j\in Y, i=0,1,...$.  
Let $A= {\rm cl\,conv}\{\dfrac{e_0}{n}+a_ne_n\,|\,n=0,1,2,\ldots\}$,
$B= {\rm cl\,conv}\{\dfrac{e_0}{n}-a_ne_n\,|\,n=0,1,2,\ldots\}$, $a_n>2^n$.
Notice that the origin belongs to the closure of $A+B$. We are going to prove that the origin does not belong to $A+B$.
If $x=\sum\limits_{j=0}^\infty x_je_j\in A$ then $x_0\geqslant 0$. Also $x_0\geqslant 0$ for all $x\in B$.
It is enough to show that $x_0> 0$ for all $x\in A$.
If $x\in {\rm conv}\{\dfrac{e_0}{n}+a_ne_n\,|\,n=0,1,2,\ldots\}$ then there exists a number $k>0$ such that $x_k\geqslant 1$.
Hence for all $x\in A$ we have $\|x\|\geqslant \dfrac{1}{C}$. Therefore, the origin does not belong to $A$. 
Suppose, $y=\sum\limits_{j=0}^\infty y_je_j\in A$ and $y_0=0$. Then $y_k>0$ for some $k$.  
There exists a sequence $(x^n)\subset {\rm conv}\{\dfrac{e_0}{n}+a_ne_n\,|\,n=0,1,2,\ldots\}$ 
such that $x^n$ tends to $y$. Then $x^n_k$ tends to $y_k$ for all $k$. We can 
represent $x^n$ as convex combination $x^n=\sum\limits_{j=0}^\infty \alpha^n_j(\dfrac{e_0}{j}+a_je_j)$.
Hence $x^n_k=\alpha^n_ka_k$ tends to $y_k$, and $\alpha^n_k$ tends to $\dfrac{y_k}{a_k}$.
Therefore, $x^n_0\geqslant \dfrac{\alpha^n_k}{k}$ where $\dfrac{\alpha^n_k}{k}$ tends to $\dfrac{y_k}{ka_k}>0$.

Notice that ${\rm cl}(A+B)= {\rm cl}({\rm conv}\{\dfrac{e_0}{n}+a_ne_n,\,|\,n=0,1,2,\ldots\}+
{\rm conv}\{\dfrac{e_0}{n}-a_ne_n\,|\,n=0,1,2,\ldots\})$.
Since for all $x\in {\rm conv}\{\dfrac{e_0}{n}+a_ne_n\,|\,n=0,1,2,\ldots\}+
{\rm conv}\{\dfrac{e_0}{n}-a_ne_n\,|\,n=0,1,2,\ldots\}$ 
we have $-a_k\leqslant x_k\leqslant a_k$, the same holds true for all $x\in {\rm cl}(A+B)$ and $k>0$.
Similarly, $0\leqslant x_0\leqslant 2$ for all $x\in {\rm cl}(A+B)$.
Therefore, if $v$ belongs to the recc\,${\rm cl}(A+B)$ then $v=0$. Hence both $A$ and $B$ have a trivial recession cone. 
\end{proof}

Notice that the sets $A$ and $B$ from the last example cannot be finitely well-positioned because in a reflexive 
Banach space no unbounded closed convex set with trivial recession cone is finitely well-positioned.

%%%%%%%%%%%%%%%%%%%%%%%%%%%%%%%%%%%%%%%%%%%%%%%%%%%%%

\section{Order cancellation law in normed spaces.}
In this section, we generalize the order law of cancellation for 
closed convex sets with common recession cone proved by Robinson in \cite{sR} 
for finite dimensional spaces.

\begin{Def}
{\rm Let $X$ be a normed space and let $A$ be a subset of $X .$ Let $\tau$ be a linear topology weaker 
than norm topology. By $\mathfrak{B}_A^{\tau}$, we define the set  
$\left\{\xi \in X\,|\,\text{there exists } (x_n )\subset A \text{ with } 
\|x_n \|\rightarrow\infty \text{ and } \dfrac{x_n }{\|x_n \|}
\rightharpoonup \xi \right\}$ where the symbol "$\rightharpoonup $" 
denotes the $\tau$-convergence in $X.$ }
\end{Def}

Obviously, the set $\mathfrak{B}_A^{\tau}$ is a subset of an asymptotic cone $A^{\tau}_{\infty}$. 
Also, the cone $A^{\tau}_{\infty}$ is a cone generated by $\mathfrak{B}_A^{\tau}$.
The definition of the set $\mathfrak{B}_A$ can be found in \cite{mmlm} for the weak topology $\tau$. 
We denote by $\mathfrak{B}_A^{\|\cdot\|}$, $\mathfrak{B}_A$, $\mathfrak{B}_A^{*}$ the set $\mathfrak{B}_A^{\tau}$ 
for, respectively, norm, weak and *-weak topology in $X$. We have $\mathfrak{B}^s_A\subset\mathfrak{B}_A
\subset\mathfrak{B}^{\ast}_A$. 
The following proposition is obvious.

\begin{Prop}\label{domkBA}
 Let $X$ be a normed space and let $A\subset X$ and $\tau$ be a linear topology weaker 
than norm topology. Then 
$ \mathfrak{B}^{\tau}_{\textup{cl} A}=\mathfrak{B}^{\tau}_A $ for the closure $\textup{cl} A$  in the norm topology.
\end{Prop}

In \cite{mmlm} it was shown that 
recc$\,A= A\dotminus A =\text{cone}\,\mathfrak{B}_A$ if the set $A$ is closed and convex. 
Now we show the following proposition.

\begin{Prop}\label{BsA}
 Let $X$ be a normed space and let $A\subset X$ be a closed and convex subset 
of $X.$ If $ \,{\rm recc}\,A\neq\{0\}$   then 
$ {\rm recc}\,A =\textup{cone}\,\mathfrak{B}^{\|\cdot\|}_A . $
\end{Prop}
\begin{proof}
Let $\xi \in {\rm recc}\,A,$ $\xi \neq 0$. 
By Lemma \ref{recconesequence} there exist sequences 
$(x_n ) \subset A $ and $(t_n)\subset \mathbb{R} ,$ $ t_n \downarrow 0$    
such that $t_n x_n \rightarrow \xi , $ hence $\dfrac{x_n }{\|x_n \|} 
=\dfrac{t_n x_n }{t_n \|x_n\|} \rightarrow \dfrac{\xi }{\|\xi \|}.$ 
Therefore,  $\dfrac{\xi }{\|\xi \|}\in \mathfrak{B}^{\|\cdot\|}_A , $ 
so  $ {\rm recc}\,A \subset \text{cone} \,\mathfrak{B}_A^{\|\cdot\|} .$
\end{proof}

\begin{Rem}
{\rm Let us notice that $\mathfrak{B}^{\|\cdot\|}_A\neq \emptyset$  implies that 
recc$\,A=A\dotminus A\neq\{0\}$. Equivalently, 
recc$\,A=\{0\}$ implies that $\mathfrak{B}^{\|\cdot\|}_A= \emptyset$}
\end{Rem}

\begin{Def}
{\rm Let $X$ be a normed space, $A\subset X$ and $\tau$ be a linear topology weaker than norm topology. 
We say that $A$ is $\tau$-narrow if for any sequence
$(x_n )\subset A \text{ with } \|x_n \|\rightarrow\infty $ we can choose a subsequence  
$(x_{n_k})$ with $\dfrac{x_{n_k}}{\|x_{n_k}\|}$ converging in the topology $\tau$ to a point other than the origin.}
\end{Def}

The following proposition is obvious implication of Proposition \ref{domkBA}.

\begin{Prop}\label{tauBA}
Let $X$ be a normed space and let $A\subset X$ and $\tau$ be a linear topology weaker 
than norm topology. Let the set $A$ be $\tau$-narrow. Then 
the closure $\textup{cl} A$  in the norm topology is also $\tau$-narrow.
\end{Prop}

\begin{proof}
Let $(x_k)\subset \textup{cl} A$ be a sequence such that $\|x_k\|$ tends to $\infty$. Then $x_k=a_k+u_k$ for some 
sequences $(a_k)\subset A$ and $(u_k)\subset X$ such that $\|u_k\|$ tends to $0$. Hence $\|a_k\|=\|x_k-u_k\|$ tends to $\infty$. 
There exists a subsequence $(a_{k_l})$ such that $\dfrac{a_{k_l}}{\|a_{k_l}\|}$ tends to some point $a\neq 0$ in the 
topology $\tau$. We have the difference $\dfrac{x_{k_l}}{\|x_{k_l}\|}-\dfrac{a_{k_l}}{\|a_{k_l}\|}
=\dfrac{u_{k_l}}{\|x_{k_l}\|}+\dfrac{a_{k_l}}{\|a_{k_l}\|}\dfrac{\|a_{k_l}\|-\|x_{k_l}\|}{\|x_{k_l}\|}$ tending to the origin. 
Therefore, $\dfrac{x_{k_l}}{\|x_{k_l}\|}$ tends to $a\neq 0$. 
\end{proof}

If a unit ball in $X$ is $\tau$-sequentially compact then a set $A$ is $\tau$-narrow if and only if $0\notin \mathfrak{B}^{\tau}_A$

\begin{tw}\label{rob01}
Let $X$ be a normed space and $\tau$ be a linear topology in $X$ weaker than the norm topology. 
Let $A,B,C\subset X$, $C$ be closed in $\tau$ and convex, $\textup{recc}C$ be pointed, 
$B$ be $\tau$-narrow and $B^{\tau}_\infty\subset \textup{recc}C$. Then the inclusion $A+B\subset B+C$ implies that $A\subset C.$
\end{tw}

\begin{proof}
Assume that $0\in B$ and $A=\{0\}$. Then by induction we can prove that for every $n\in\mathbb{N}$ we have 
$0\in B\subset B+\underbrace{C+...+C}_{n}= B+nC$ hence $0\in \bigcap\limits_{n\in\mathbb{N}} \Big(\dfrac{B}{n} +C\Big).$
Then there exist $b_n \in B$ and $c_n\in C$ such that $0=\dfrac{b_n}{n} +c_n .$

Case 1. The sequence $\dfrac{b_n}{n}$ is norm-bounded. If $(b_n)$ is norm-bounded  than $\dfrac{b_n}{n}$ tends to 
$0$ and $0$ is the limit of the sequence $(c_n)$. If $(b_n)$ is not norm-bounded then,
applying the assumption that the set $B$ is $\tau$-narrow, for some subsequence $(b_{n_k})$
the sequence $\dfrac{b_{n_k}}{\|b_{n_k}\|}$ tends in the topology $\tau$ to some $b_0\in B^{\tau}_{\infty}\setminus \{0\}$. 
The number sequence $\dfrac{\|b_{n_k}\|}{{n_k}}$ is bounded and contains a subsequence convergent to $t\geqslant 0$. 
Hence some subsequence $\dfrac{b_{n_{k_l}}}{n_{k_l}}=\dfrac{b_{n_{k_l}}}{\|b_{n_{k_l}}\|}\dfrac{\|b_{n_{k_l}}\|}{n_{k_l}}$ 
tends in the topology $\tau$ to $tb_0, t\geqslant 0$. 
Then $\lim_{l\rightarrow \infty}c_{n_{k_l}}  = -tb_0 \in C.$
Notice that $tb_0\in\,B^{\tau}_\infty\subset\textup{ recc}\,C $, and $0=-tb_0 +tb_0 \in C.$ 

Case 2. The sequence $\dfrac{b_n}{n}$ is not bounded. 
We may assume that  $\|\dfrac{b_n}{n}\|\longrightarrow\infty$
and write 
$0=\dfrac{b_n}{\|b_n\|} +\dfrac{nc_n }{\|b_n\|} .$
By assumption there exists a subsequence 
$\dfrac{b_{n_k}}{\|b_{n_k} \|}$ tending in the topology $\tau$ to some $b_0 \neq 0$ from $B^{\tau}_\infty$
and, consequently, $\lim_{k\rightarrow \infty}\dfrac{n_k c_{n_k} }{\|b_{n_k}\|}=-b_0 .$
We obtain $b_0  \in B^{\tau}_\infty \subset {\rm recc}\,C$ and, by  Lemma \ref{recconesequence}, also 
$-b_0  \in C^{\tau}_\infty={\rm recc}\,C$ which is impossible, since recc\,{C} is pointed.

We just have proved that the theorem holds in the case of $0\in B$ and $A=\{0\}$. In general case 
let $A+B\subset B+C$ and $b\in B$. Consider any $x\in A$. Then $x+(B-b)\subset (B-b)+C$, and $\{0\}+(B-b)\subset (B-b)+(C-x)$. 
The sets $\{0\}$, $(B-b)$ and $(C-x)$ satisfy assumptions of the theorem in its proved case.
Hence $0\in (C-x)$, and $x\in C$. Therefore, $A \subset C$.   
\end{proof}

\begin{Lem}\label{rob01lem}
Let $X$ be a locally convex space. Let for $B\subset X$ and any
$A,C\subset X$, such that $C$ be closed and convex,  
the inclusion $A+B\subset B+C$ imply that $A\subset C.$
Then the inclusion $A+B\subset \textup{cl}(B+C)$ imply that $A\subset C.$ 
\end{Lem}

\begin{proof}
Notice that the implication $A+B\subset B+$cl$(C+U)\Longrightarrow A\subset $cl$(C+U)$ holds true 
for any convex neighborhood $U$ of the origin. Assume that $A+B\subset \textup{cl}(B+C)$.
We obtain that $A+B\subset B+C+U\subset B+$cl$(C+U)$ and, by the assumptions of the theorem, $A\subset $cl$(C+U)\subset C+U+U.$   
Then $A+B\subset \textup{cl}(B+C)$ implies that $A$ is contained in the intersection of all $C+2U$. Therefore, $A\subset C$.
\end{proof}

\begin{Cor}\label{rob01cor}
Let $X$ be a normed space. Assume that there exists a locally convex topology $\tau$, weaker than the norm topology.
Let $A,B,C\subset X$, $C$ be closed in $\tau$ and convex, $B$ be $\tau$-narrow, $\textup{recc}C$ be pointed and 
$B^{\tau}_\infty\subset \textup{recc}C$. 
Then the inclusion $A+B\subset \textup{cl}_{\tau}(B+C)$ implies that $A\subset C.$
\end{Cor}

\vspace{2mm}
\noindent
A normed space $X$ is said to be a {\sc wcf} space (see \cite{AL}) if 
it contains a weakly compact subset $K$, such that $X={\rm cl\,lin}K$.  

\begin{Cor}\label{rob1cor2}
Let $A,B,C$ be subsets of a dual $X^{\ast}$ to a 
{\sc wcf} space $X$, $C$ be *-weakly closed and convex, $B$ be nonempty and $B_\infty^*$ be asymptotic cone with respect 
to *-weak topology. Assume that $0\notin \mathfrak{B}^{\ast}_B $,  
and $B_\infty^*\subset \text{\rm recc}\, C$. 
If the cone $\text{\rm recc}\, C$ is pointed then the inclusion $A+B\subset \textup{cl}^*(B+C)$ implies that $A\subset C.$
\end{Cor}
\begin{proof}
Since $X$ is {\sc wcf}, by Corollary 2 in \cite{AL} and Eberlein--Smulian theorem \cite{rW} 
the unit ball in $X^{\ast}$ is *-weakly sequentially compact, and $B$ is *-weakly-narrow.
Applying Corollary \ref{rob01cor} we obtain our corollary.
\end{proof}

Since every reflexive Banach space $X$ is a dual of a {\sc wcf} space and its *-weak 
topology coincides with the weak topology, 
the next theorem follows from Corollary \ref{rob1cor2}.

\begin{tw}\label{rob2}
Let $A,B,C$ be subsets of a reflexive Banach space $X$, $C$ be closed and convex, $B$ be nonempty 
and $B_\infty$ be asymptotic cone with respect 
to weak topology. Assume that $0\notin \mathfrak{B}_B $, 
$B_\infty \subset\text{\rm recc}\,C =V$ and $V$ is a pointed cone.  
Then the inclusion $A+B\subset \textup{cl}_{\textup{weak}}(B+C)$ implies that $A\subset C.$
\end{tw}

The theorem generalizes Robinson's result (Lemma 1 in \cite{sR}) 
in $\mathbb{R}^n$. 
The condition $0\notin \mathfrak{B}_B $ in the theorem is crucial. 
Without it the theorem is false in every infinite dimensional Banach space. In Theorem
\ref{reccA+B} we have recc$\,(A\,\dot{+}\,B)=[0,\infty)e_0$. Then $\big([0,\infty)e_0+B\big)+A\subset A\,\dot{+}\,B$. 
Yet, $\big([0,\infty)e_0+B\big)\not\subset B$.
The reason is that  
even though recc$\,A=$ recc$\,B=\{0\}$ and $\mathfrak{B}_A^s=\emptyset$, 
we have $0\in \mathfrak{B}_A $.  

\begin{Prop}\label{clB+C}
Let $X$ be a reflexive Banach space. 
Let $B,C\subset X$, $B,C$ be closed and convex, $(\textup{recc}B)\cap(-\textup{recc}C)=\{0\}$ 
and $0\notin \mathfrak{B}_B $\textup{(}$B$ be weakly-narrow\textup{)}. 
Then the sum $B+C$ is closed. 
\end{Prop}

\begin{proof}
Let $x\in \textup{cl}(B+C)$. Then there exist sequences 
$(b_n)\subset B$, $(c_n)\subset C$ and $(u_n)\subset X, \|u_n\|\downarrow 0$ such that 
$b_n+c_n+u_n=x$.
First, if $\|b_n\|$ is bounded then some subsequence $(b_{n_k})$ weakly tends to $b\in B$. 
Hence $c_{n_k}$ tends to $x-b\in C$, and $x\in B+C$.
Second, if $\|b_n\|$ tends to infinity then 
some subsequence $\dfrac{b_{n_k}}{\|b_{n_k}\|}$ weakly tends to some $b\in$  recc\,$B, b\neq 0$. 
Hence $\dfrac{c_{n_k}}{\|c_{n_k}\|}$ tends to some $-b\in$  recc\,$C$, which contradicts the assumption 
that the cone $V$ is pointed. 
\end{proof}

Proposition \ref{clB+C} is also a consequence of Dieudonn\'e's theorem \cite{JDi}.

\begin{Prop}\label{reccB+.+C}
Let $X$ be a normed vector space. 
Let $B,C\subset X$, $B,C$ be closed and convex, $\textup{recc}B=\textup{recc}C$ be pointed, 
and $B$ be weakly-narrow. 
Then $\textup{recc cl}(B+C)=\textup{recc}C$. 
\end{Prop}

\begin{proof}
Let $x\in \textup{recc cl}(B+C), x\neq 0$. Then $b_0+c_0+[0,\infty)x\subset \textup{cl}(B+C)$ for 
some $b_0\in B, c_0\in C$. There exist sequences 
$(b_n)\subset B$, $(c_n)\subset C$ and $(u_n)\subset X, \|u_n\|\downarrow 0$ such that 
$b_n+c_n+u_n=b_0+c_0+nx$.
Notice that $$x=\dfrac{b_n}{n}+\dfrac{c_n}{n}+\dfrac{u_n}{n}-\dfrac{b_0}{n}-\dfrac{c_0}{n}.$$
First, if $\dfrac{\|b_n\|}{n}$ tends to 0 then $\dfrac{c_n}{n}$ tends to $x$, and $x\in C_{\infty}=$ recc\,$C$.
Second, if $\dfrac{\|b_n\|}{n}$ tends to $\beta\in (0,\infty)$ then 
some subsequence $\dfrac{b_{n_k}}{n_k}$ tends to $\beta b\in$  recc\,$B$, where $\dfrac{b_{n_k}}{\|b_{n_k}\|}$ 
weakly tends to some point $b\neq 0$. Hence
$\dfrac{c_{n_k}}{n_k}$ tends to $x-\beta b\in$ recc\,$C$, and $x\in $ recc\,$B+$ recc\,$C=$ recc\,$C$.

Third, if $\dfrac{\|b_n\|}{n}$ tends to infinity then again some subsequence $\dfrac{b_{n_k}}{\|b_{n_k}\|}$ 
weakly tends to a nonzero element $b\in$ recc\,$B$.
Notice that also $\dfrac{b_0+c_0-u_{n_k}+n_kx-b_{n_k}}{\|b_{n_k}\|}$ tends to $-b$ while the fraction
$\dfrac{\|b_0+c_0-u_{n_k}+n_kx-b_{n_k}\|}{\|b_{n_k}\|}$ tends to $1$.
Then $\dfrac{c_{n_k}}{\|c_{n_k}\|}=\dfrac{b_0+c_0-u_{n_k}+n_kx-b_{n_k}}{\|b_0+c_0-u_{n_k}+n_kx-b_{n_k}\|}$ 
tends to $-b$ and $-b\in$ recc\,$C$, which contradicts the assumption that the cone recc\,$C$ is pointed.
\end{proof}

Obviously, the assumption of weak-narrowness of one of the sets in the last proposition is essential.

\begin{tw}\label{B+.+C}
Let $X$ be a normed space and $\tau$ be a linear topology weaker than the norm topology.
Let subsets $B,C\subset X$ be $\tau$-narrow and $B^{\tau}_{\infty}\cap (-C^{\tau}_{\infty})=\{0\}$. 
Then the set $\textup{cl}_{\|\cdot\|}(B+C)$ is $\tau$-narrow. 
\end{tw}

\begin{proof}
Let a $(x_k)\subset B+C$,  $\|x_k\|$ tends to $\infty$. For some sequences $(b_k)\subset B$ and $(c_k)\subset C$ we have $x_k=b_k+c_k$.

First, assume that the sequence $(\|b_k\|)$ is bounded. Hence $\|c_k\|\longrightarrow \infty$. 
For some subsequence $(c_{k_l})$ terms $\dfrac{c_{k_l}}{\|c_{k_l}\|}$ tend to some $c\neq 0$ in the topology $\tau$.  
Also $\dfrac{x_{k_l}}{\|x_{k_l}\|}-\dfrac{c_{k_l}}{\|c_{k_l}\|}=
\dfrac{x_{k_l}}{\|x_{k_l}\|}-\dfrac{c_{k_l}}{\|x_{k_l}\|}+\dfrac{c_{k_l}}{\|x_{k_l}\|}-\dfrac{c_{k_l}}{\|c_{k_l}\|}$  
$=\dfrac{b_{k_l}}{\|x_{k_l}\|}+\dfrac{\|c_{k_l}\|-\|x_{k_l}\|}{\|x_{k_l}\|}\dfrac{c_{k_l}}{\|c_{k_l}\|}$. 
Since the sequence $(\|x_{k_l}-c_{k_l}\|=\|b_{k_l}\|)$ is bounded, 
a term $\dfrac{\|c_{k_l}\|-\|x_{k_l}\|}{\|x_{k_l}\|}$ tends to 0. Since $\dfrac{c_{k_l}}{\|c_{k_l}\|}$ is convergent 
in the topology $\tau$, the sequence 
$\dfrac{x_{k_l}}{\|x_{k_l}\|}-\dfrac{c_{k_l}}{\|c_{k_l}\|}$ tends to 0. Then the sequence $\dfrac{x_{k_l}}{\|x_{k_l}\|}$ 
tends to $c\neq 0$.  

Second, let $\|b_k\|\longrightarrow \infty$ and $\|c_k\|\longrightarrow \infty$. 
By assumption there exists a sequence $(k_l)$ such that $\dfrac{b_{k_l}}{\|b_{k_l}\|}\longrightarrow b\neq 0$, 
$\dfrac{c_{k_l}}{\|c_{k_l}\|}\longrightarrow c\neq 0$ and $\dfrac{\|b_{k_l}\|}{\|b_{k_l}\|+\|c_{k_l}\|}\longrightarrow t\in [0,1]$.
If $\dfrac{\|b_{k_l}\|+\|c_{k_l}\|}{\|x_{k_l}\|}\longrightarrow  \infty$ 
then $\dfrac{x_{k_l}}{\|b_{k_l}\|+\|c_{k_l}\|}=\dfrac{b_{k_l}}{\|b_{k_l}\|}\dfrac{\|b_{k_l}\|}{\|b_{k_l}\|+\|c_{k_l}\|}+
\dfrac{c_{k_l}}{\|c_{k_l}\|}\dfrac{\|c_{k_l}\|}{\|b_{k_l}\|+\|c_{k_l}\|}$
is convergent and tends to $0=tb+(1-t)c$. Then $b=0$ or $c=0$ or $0\neq tb=(t-1)c\in B^{\tau}_{\infty}\cap(-C^{\tau}_{\infty})$ 
which contradicts the assumptions.

Therefore, we have $\dfrac{\|b_{k_l}\|+\|c_{k_l}\|}{\|x_{k_l}\|}\longrightarrow r\in [1,\infty)$.
Then the term $\dfrac{x_{k_l}}{\|x_{k_l}\|}=\left(\dfrac{b_{k_l}}{\|b_{k_l}\|}\dfrac{\|b_{k_l}\|}{\|b_{k_l}\|+\|c_{k_l}\|}+
\dfrac{c_{k_l}}{\|c_{k_l}\|}\dfrac{\|c_{k_l}\|}{\|b_{k_l}\|+\|c_{k_l}\|}\right)\dfrac{\|b_{k_l}\|+\|c_{k_l}\|}{\|x_{k_l}\|}$
tends to $r(tb+(1-t)c)\neq 0$.

We have just proved that $B+C$ is $\tau$-narrow. By Proposition \ref{tauBA} the set cl$_{\|\cdot\|}(B+C)$ is $\tau$-narrow.
\end{proof}

\begin{Cor}\label{B+.+Crefl}
Let $X$ be a reflexive Banach space. Let $B,C\subset X$ be closed and convex sets sharing a pointed recession cone 
$\textup{recc}B=\textup{recc}C$. Assume that $0\notin \mathfrak{B}_B\cup \mathfrak{B}_C$. 
Then $0\notin \mathfrak{B}_{\textup{cl}(B+C)}$, where $\textup{cl}$ is the weak closure. 
\end{Cor}

\begin{proof}
Notice that for reflexive spaces the unit ball in $X$ is weakly sequentially compact. Hence the condition 
$0\notin \mathfrak{B}_B\cup \mathfrak{B}_C$ implies that $B$ and $C$ are weakly narrow. Since their asymptotic cones 
are contained in one pointed recession cone, the assumptions of Theorem \ref{B+.+C} are fulfilled for 
the weak topology $\tau$. 
Therefore, the set $\textup{cl}_{\|\cdot\|}(B+C)=\textup{cl}(B+C)$ is weakly narrow and $0\notin \mathfrak{B}_{\textup{cl}(B+C)}$.
\end{proof}

\begin{tw}\label{mrh1}
Let $X$ be a normed vector space and $V\subset X$ be a closed convex pointed cone. Then the family of all weakly narrow 
closed convex subsets of $X$ sharing the recession cone $V$ with modified Minkowski addition $\dot{+}$ defined by 
$A\dot{+} B:=\textup{cl}\,(A+B)$ and multiplication by non-negative numbers $($we assume $0\cdot A:=V)$ 
is an abstract convex cone having a neutral element $V$ and satisfying cancellation law. 
\end{tw}

\begin{proof}
First, notice that for any two weakly narrow closed convex sets $A$ and $B$ with a common pointed recession cone $V$ 
the set $\textup{cl}\,(A+B)$ is weakly narrow by Theorem \ref{B+.+C} and the recession cone of $\textup{cl}\,(A+B)$ 
is equal to $V$ by Proposition \ref{reccB+.+C}. Then the family of all weakly narrow 
closed convex sets is closed with respect to the addition $\dot{+}$. Checking conditions for an abstract convex cone 
is straightforward. The cancellation law in the family of all weakly narrow 
closed convex sets follows from Theorem \ref{rob2}.
\end{proof}

\noindent
In view of Proposition \ref{clB+C} we can formulate the following Corollary. 

\begin{Cor}\label{mrh1+}
Let $X$ be a reflexive Banach space and $V$ be a closed convex pointed cone. Then the family of all subsets $A\subset X$ sharing the  
recession cone $V$ and satisfying the condition $0\notin \mathfrak{B}_A$ with Minkowski addition 
and multiplication by non-negative numbers $($we assume $0\cdot A:=V)$ 
is an abstract convex cone having a neutral element $V$ and satisfying cancellation law. 
\end{Cor}

\noindent
Let $X$ be a normed vector space and $V\subset X$ be a closed convex pointed cone. Then by Theorem \ref{mrh1} 
we can embed the family $\mathcal{C}_V^{\textup{n}}(X)$ of all weakly narrow 
closed convex subsets of $X$ sharing the recession cone $V$ in a Minkowski--R\aa dstr\"om--H\"ormander vector 
space $(\mathcal{C}_V^{\textup{n}}(X))^2/_{\sim}$. Moreover, by Theorem 2.5 in \cite{GP} if, additionally, 
$X$ is a reflexive Banach space and the set $\{x\in V\,|\,\|x\|=1\}$ can be strictly separated from the origin 
then every quotient class $[A,B]\in(\mathcal{C}_V^{\textup{n}}(X))^2/_{\sim}$ has a minimal representative.

%\\

%\\
%\noindent
%Let $\mathbb{B}$ be a closed unit ball in a reflexive Banach space $X$ and $\mathbb{S}$ be a unit sphere. Denote
%$V_n:=$cl\,cone\,(conv$(V\cap \mathbb{S})+\dfrac{1}{n}\mathbb{B})$. If $V$ is a pointed closed convex cone 
%($V\cap \mathbb{S}$ can be separated from the origin by closed hyperplane $H$ ($V_n\subset $cone($\mathbb{B}\cap tH$) )
%(find example clconeconv$\{e_n+e_1/n,-e_n+e_1/n\}$)?) then for some $n_0$ all cones $V_n, n\geqslant n_0$ are pointed.

%\noindent
%Let $d_n(A,B):=d_H(A+V_n,B+V_n)$ for $A,B\in \mathcal{C}_V(X)$.
%$d(A,B)=\sum\limits_{n=n_0}^{\infty}2^{-n}\dfrac{d_n(A,B)}{1+d_n(A,B)}$

\section{Order cancellation law in topological vector spaces}
In this section we prove order cancellation law for topological vector space where we cancel by closed and convex sets . 
These results generalize theorems obtained by  and Tabor, Bielawski in \cite{BT}.

Let $X$ be a topological vector space and let $V\subset X$ be a closed 
and convex cone.
By $\mathcal{C}(X)$ we denote the family of al nonempty, closed 
and convex subsets of $X$ 
and by $\mathcal{B}(X) $ the family of all nonempty, bounded, closed and convex  subsets of $X.$ 
Now let 
$\mathcal{C}_V^0 (X) =\{A\in \mathcal{C} (X) \,|\, V\subset \mbox{recc}\,A \}  $
and
\begin{equation}
\mathcal{C}_V (X) =\{A\in \mathcal{C} (X) \,|\, A\subset V\,
\dot{+}\,B,\,\, V\subset A\,\dot{+}\,B \mbox { for some } 
B\in \mathcal{B} (X) \}.  \tag{5.1}
\end{equation}

We may look at elements of $\mathcal{C}_V (X)$ as at sets of bounded Hausdorff-like distance from the cone $V$. 
The family $\mathcal{C}_V (X)$ is a subfamily  of all 'bounded' elements of $\mathcal{C}_V^0 (X)$. 
For a sequence $(A_n )\subset \mathcal{C}_V^0 (X)$ we define the limit operator as follows: 
$\lim(A_n ) :=A\in \mathcal{C}_V^0 (X)$ if and only if for every neighborhood $U$ of zero in $X$ 
there exists $k\in \mathbb{N}$ such that $A\subset A_n + U \mbox{ and } A_n\subset A +U$ for $n\geqslant k.$ 

Moreover, the family $\mathcal{C}_V^0 (X)$ with the addition 
$A\,\dot{+}\,B=\textup{cl}(A+B) $ and the multiplication by 
non-negative numbers defined by $\lambda A:=\{\lambda a \,|\, a\in A\},$ 
$0\cdot A :=V$ is an abstract convex cone. 
\begin{Prop}
The algebraic system $(\mathcal{C}_V^0(X), \cdot \,, \dot{+} , \subset , \lim )$ satisfies the axioms {\rm (S1--11)}.
\end{Prop}
\begin{proof} Obviously, the axioms (S1--2) and (S4--9) 
are satisfied. So we need to prove only the axioms (S3) and (S10--11). 
Since $V$ is a neutral element in $\mathcal{C}_V^0 (X)$, the axiom (S--3) holds true.

To prove (S10) let us assume that $\lim (A_n)=A$, $\lim(B_n )=B$ and $A_n\subset B_n $ for all $n\in\mathbb{N} .$ 
Take any neighborhood $U$ of the origin in $X$ then $A\subset A_k +U \subset B_k +U \subset B+U +U $ 
for sufficiently large $k$.
Hence $$A\subset \bigcap\{(B+U+U )\,|\,U\textup{ is a neighborhood of } 0\} =\textup{cl}B =B,$$ so the axiom (S10) is satisfied.

To prove (S11), let us assume that $\lim(A_n )=A$, $\lim(B_n )=B.$ Take any neighbourhood $U$ of zero in $X$ 
and let $W$ be a neighborhood of zero in $X$ such that $W+W+W\subset U.$ Let $k$ be a natural number such that 
$A\subset A_n + W$, $A_n\subset A +W$, $B\subset B_n + W$ and $B_n\subset B +W$ 
for  $n\geqslant k.$ 
Then $$A\,\dot{+}\,B\subset A+B+W\subset A_n + B_n  + W +W+W 
\subset A_n +B_n +U\subset A_n\, \dot{+}\,B_n +U$$
and, analogously, 
$ A_n \dot{+}B_n\subset A\dot{+}B+U$ for  
$n\geqslant k$. Hence $\lim (A_n \dot{+}B_n )= A\dot{+}B .$ Therefore, the axiom (S11)  is satisfied. 
\end{proof}
Now, we prove two lemmas concernig the family $\mathcal{C}_V (X) $
\begin{Lem}
For any $A\in \mathcal{C}_V (X) $ we have {\rm ${ \mbox{recc}}A =V.$}
\end{Lem}
\begin{proof}
Take any $A\in \mathcal{C}_V (X)$. First, we prove that $V\subset \mbox{recc} A$.  
Let us observe that for every $n\in\mathbb{N}$ we have $V= \dfrac{1}{n} V\subset\dfrac{1}{n}A \,\dot{+}\,\dfrac{1}{n} B$
for some $ B\in \mathcal{B} (X).$ 
Now, take arbitrary $a\in A , v\in V$. For any neighborhood $U$ of the origin let $W$ be another neighborhood of the origin 
such that $W+W+W \subset U.$ Let also $m\in \mathbb{N}$ be large enough that $\dfrac{a}{m}\in W$ and $\dfrac{1}{m} B\subset W.$ 
Then for some $a_1\in A $ and $b_1\in B$ we have 
\begin{multline*}
%begin{split}
a+v\in  a+\dfrac{1}{m} a_1 +\dfrac{1}{m} b_1 + W=\left[\dfrac{1}{m} a_1 
+\left(1-\dfrac{1}{m}\right) a\right]+\dfrac{1}{m} a +\dfrac{1}{m} b_1 +W\,\,\,\,\,\,\\  
\subset A+W+W+W\subset A+U.\hspace{3cm}
%\end{split}
\end{multline*}
Hence $A+V\subset \textup{cl}A =A.$ Then $V\subset \mbox{recc}\,A .$ 

To prove the inverse inclusion observe that $A\subset V\,\dot{+}\,B $ for some $B\in\mathcal{B} (X) $. Take any $a\in A$. We have 
$$\mbox{recc}\,A=\mbox{recc}\,(A-a) 
=\bigcap_n \dfrac{A-a}{n}\subset \bigcap_n \dfrac{V\dot{+}B-a}{n}
=\bigcap_n (V\dot{+} \dfrac{B-a}{n}) \subset \textup{cl}V =V,$$ 
and the proof is complete.
\end{proof}
\begin{Lem}\label{ogrcv}
For any $A\in \mathcal{C}_V (X) $ we have $\lim(2^{-n} A)=V.$

\end{Lem}
\begin{proof}
Let $B\in\mathcal{B} (X) $ be such that $A\subset V\,\dot{+}\, B,\,\, V\subset A\,\dot{+}\, B.$ Take any neighborhood $U$ of zero in $X$ 
and let $W$ be a balanced neighborhood of zero such that $W+W\subset U$ and let $k\in\mathbb{N} $ be such that $2^{-n} B\subset W $ 
for $n\geqslant k .$ Then we have $2^{-n } A \subset 2^{-n} ( V+ B +W ) \subset V +W +W \subset V +U $ and, analogously, 
$V=2^{-n}V\subset 2^{-n } A +U .$ Thus $\lim(2^{-n} A)=V.$
\end{proof}

The following theorems follow from the above considerations and Section 2.

\begin{tw}\label{gTB0}
Let $A,B, C\in \mathcal{C}_V^0 (X)$ and suppose that $\lim(2^{-n}B)$ exists and is contained in $V$ 
then the inclusion $A\,\dot{+}\, B \subset B\,\dot{+}\, C $ implies $A\subset C.$
\end{tw}

\begin{proof}
In order to apply Lemma \ref{olc<<} notice that convexity of a closed convex set $C$ implies the convexity in the sense of 
Definition \ref{convexity}.
\end{proof}

The next theorem is a straightforward corollary from Theorem \ref{gTB0} and Lemma \ref{ogrcv}.  

\begin{tw}\label{gTB}
Let $A,B, C\in \mathcal{C}_V (X)$  then the inclusion $A\,\dot{+}\, B \subset B\,\dot{+}\, C $ implies $A\subset C.$
\end{tw}
Tabor and Bielawski \cite{BT} proved a version of order cancellation law for closed convex sets having finite Hausdorff 
distance from a fixed convex cone $V$ in a normed space $X$.
Theorem \ref{gTB} is another version of order cancellation law which generalizes to topological spaces 
a result obtained by Tabor and Bielawski.  
 
In the Cartesian product $ \mathcal{C}_V (X)\times \mathcal{C}_V (X)$, due to cancellation law 
we can introduce an equivalence relation $\sim $ in an usual way:
\begin{align*}
(A,B) \sim (C,D) \Longleftrightarrow A\,\dot{+}\, D =B\,\dot{+}\, C
\end{align*}
Let $\widetilde{\mathcal{C}}_V (X) =\mathcal{C}_V (X)\times \mathcal{C}_V (X)/_{\sim}$ be the quotient set of the equivalence classes
with respect to $\sim\ $. By $[A,B]\in \widetilde{\mathcal{C}}_V (X)$ we denote the equivalence class of an element 
$(A,B)\in  \mathcal{C}_V (X)\times \mathcal{C}_V (X).$

The semigroup $\mathcal{C}_V (X)$ with the addition $\dot{+}$ and scalar multiplication defined by $\lambda A:=\{\lambda a:a\in A\}$ for 
$\lambda >0$ and $0\cdot A=V$ is an abstract convex cone. 
Thanks to Theorem \ref{gTB},
it can be proved that the quotient set $ \widetilde{\mathcal{C}}_V (X)$ with addition 
\begin{align}
[A,B] +[C,D] :=[A\,\dot{+}\, C , B\,\dot{+}\, D]
\end{align}
and scalar multiplication defined by 
\begin{align}
\lambda [A, B] :=\begin{cases} [\lambda A , \lambda B] \mbox{ if } \lambda \geqslant 0\\ 
[(-\lambda )B, (-\lambda )A ] \mbox{ if } \lambda \leqslant 0\end{cases}
\end{align}
is a real vector space.

Moreover, the function $j: \mathcal{C}_V (X) \to \widetilde{\mathcal{C}}_V (X)$ defined by $j (A) = [A, V]$ is 
a canonical embedding which satisfies $j(A\,\dot{+}\, B ) =j(A) +j(B) \mbox{ and } j(\lambda A) =\lambda j(A) $
for all $A,B\in  \mathcal{C}_V (X)$ and $\lambda \geqslant 0.$

Now, let $\mathcal{W}$ be a basis of balanced neighbourhoods of $0$ in topological vector space $X.$ 
For $U \in \mathcal{W}$ we define by $W_U$ the set  
$$\{ (A,B) \in (\mathcal{C}_V (X))^2 : A\subset B\,\dot{+}\,K , 
B\subset A\,\dot{+}\,K \mbox{ for some } 
K \in \mathcal{C} (X) , K\subset U\}.$$

The family $\{W_U \}_{U\in\mathcal{W}}\,$ forms a basis of uniformity 
on $\mathcal{C}_V (X)$. The uniformity induces a topology in $\mathcal{C}_V (X) $. Let us denote this topology by $\tau_H$. 
The scalar multiplication $\cdot$ and the addition $\dot{+} $ in $\mathcal{C}_V (X)$ are continuous in the topology $\tau_H$.

We define a basis of neighbourhoods of $0$ in the space $ \widetilde{\mathcal{C}}_V (X)$ as follows. 
Let $U\in \mathcal{W} $ and let $$\mathcal{T}_U 
=\{ [A,B]\in\widetilde{\mathcal{C}}_V (X):
A\subset B\,\dot{+}\,K , B\subset A\,\dot{+}\,K \mbox{ for some } 
K \in \mathcal{C} (X) , K\subset U\}.$$ 
The family  $\{\mathcal{T}_U \}_{U\in\mathcal{W}}\,\,$ forms a basis of neighbourhoods of $0$ in  the space $ \widetilde{\mathcal{C}}_V (X)$ 
and defines a linear topology on the space  the space 
$ \widetilde{\mathcal{C}}_V (X)$. Denote this topology by $\sigma_H .$

Now, we can express the following theorem which generalizes embedding theorem obtained by Tabor and Bielawski  for normed spaces \cite{BT}.

\begin{tw}\label{BT2}
The canonical mapping $ j: \mathcal{C}_V (X) \to j(\mathcal{C}_V (X))
\subset \widetilde{\mathcal{C}}_V (X)$ is an isomorphic embedding 
of a topological abstract 
convex cone $(\mathcal{C}_V (X) , \dot{+} , \cdot , \tau_H ) $ into 
a real topological vector space $(\widetilde{\mathcal{C}}_V (X) ,+,\cdot 
,\sigma_H )$. Moreover the canonical mapping $ j: \mathcal{C}_V (X) \to 
j(\mathcal{C}_V (X))\subset \widetilde{\mathcal{C}}_V (X)$ is a homeomorphism.
\end{tw} 

The results from this section apply also to topological vector spaces like $l^p, L^p, 0<p<1$. In the case of 
locally convex vector spaces the introduced topology in $\mathcal{C}_V(X)$ coincides with "Hausdorff" topology, where 
sets of families $\mathcal{T}_{A,U}:=\{B\in \mathcal{C}_V(X)\,|\,A\subset B\,\dot{+}\,U,B\subset A\,\dot{+}\,U\},U\in \mathcal{W}$ 
form basis of neighborhoods of $A\in \mathcal{C}_V(X)$. In spaces which are not locally convex a problem arises 
because belonging $B\in \mathcal{T}_{A,U}$ does not imply that $B_1\in \mathcal{T}_{A_1,U}$ for all $(A_1,B_1)\in [A,B]$.
Notice that $A\subset B\,\dot{+}\,U\Rightarrow A\,\dot{+}\,C\subset B\,\dot{+}\,C\,\dot{+}\,U,C\in \mathcal{C}_V(X)$ 
but having the inclusion $A\,\dot{+}\,C\subset B\,\dot{+}\,C\,\dot{+}\,U$
we may not be able to cancel $C$, since the set $B\,\dot{+}\,C\,\dot{+}\,U$ may not be convex. 
Then the Minkowski--R\aa dstr\"om--H\"ormander space $\mathcal{C}_V(X)\times\mathcal{C}_V(X)/_{\sim}$ 
does not posses a topology compatible with the  
topology of neighborhoods $\mathcal{T}_{A,U}$.

\section{Conclusions}
The order law of cancellation $A+B\subset B+C\Longrightarrow A\subset C$ for subsets of a vector space is a valuable tool 
in mathematical analysis. However,
proving the order cancellation law 
for an unbounded set $B$ poses a difficult challenge. 
Assuming the inclusion of recession cones recc$\,B\subset$ recc$\,C$ and the convexity of $C$ is obvious. 
In section 5 we proved that certain kind of "finite distance" defined in (5.1) between a convex set  
$B$ and recc$\,B\subset $ is sufficient. 
Moreover, the semigroup of closed convex sets sharing their recession cones can be topologically embedded 
in the quotient topological vector space. 

In the case of "infinite distance"
between $B$ and recc$\,B$ we were able to prove the order cancellation law   
under conditions that (1) the space $X$ is normed and 
(2) the set $B$ is narrow according to Definition 4.5, which means that the asymptotic cone of $B$ is generated 
by a subset of the unit sphere which can separated from the origin.   
We expect that the semigroup of closed convex narrow sets sharing their recession cones is metrizable and can be topologically embedded 
in the quotient linear metric space.

\noindent
Jerzy Grzybowski\\
Adam Mickiewicz University, Faculty of Mathematics and Computer Science, Uniwersytetu Pozna\'nskiego 4, 61-614  Pozna\'n, Poland \\
\textit{E-mail address}: jgrz@amu.edu.pl\\
\\

\noindent
Hubert Przybycie\'n\\
Jan Kochanowski University in Kielce, Department of Mathematics, Uniwersytecka 7, 25-406  Kielce, Poland \\
\textit{E-mail address}: hprzybycien@ujk.edu.pl

\end{document}